# QUESTIONNEMENTS AUTOUR DE LA SYNCHRONISATION DANS L'ENSEIGNEMENT DES MATHEMATIQUES A DES ELEVES SOURDS

ASSUDE[*] Teresa – MILLON-FAURE[**] Karine – TAMBONE[***]Jeannette

**Résumé** – Dans cette communication, nous étudions un dispositif mis en place afin de scolariser des élèves sourds et nous analysons plus particulièrement plusieurs séances de mathématiques auxquelles ont participé ces élèves. Nous montrons notamment que si tous les élèves semblent globalement en phase, des cycles de désynchronisations/resynchronisations se succèdent grâce notamment à l'apparition de bulles de compréhension.

**Mots-clefs** : élèves sourds, mathématiques, synchronisation, système didactique principal et auxiliaire.

**Abstract** – In this communication we study a device set up to school deaf pupils. We analyze some sessions of mathematics classroom in which participated these pupils. We show in particular that if all the pupils seem globally in phase, cycles of desynchronizations / resynchronizations follow one another due to the appearance of « bubbles of understanding ».

## I.  INTRODUCTION

En France, la loi n°2005-102 du 11 février 2005 *pour l'égalité des droits et des chances, la participation et la citoyenneté des personnes handicapées* dispose que les élèves en situation de handicap doivent être prioritairement scolarisés dans le milieu scolaire ordinaire. Cette volonté a été confirmée dans la loi n°2013-595 du 8 juillet *d'orientation et de programmation pour la refondation de l'école de la République*. Dans ce sens, des ULIS (Unité localisée d'inclusion scolaire) sont des dispositifs de scolarisation des élèves en situation de handicap en tant qu'une des modalités de l'accessibilité pédagogique : « Les élèves orientés en ULIS sont ceux qui, en plus des aménagements et adaptations pédagogiques et des mesures de compensation mis en œuvre par les équipes éducatives, nécessitent un enseignement adapté dans le cadre de regroupements ». La circulaire n°2015-129 du 21 août 2015 précise les modalités et l'organisation de ces dispositifs : « Les ULIS constituent un dispositif qui offre aux élèves qui en bénéficient une organisation pédagogique adaptée à leurs besoins ainsi que des enseignements adaptés dans le cadre de regroupements et permet la mise en œuvre de leurs projets personnalisés de scolarisation ».

A la rentrée de septembre 2016, une ULIS-Collège destinée à des élèves sourds et malentendants a ouvert dans un collège de Marseille. Le parcours qui est proposé aux élèves de cette ULIS est un parcours bilingue, français écrit et langue des signes française (LSF). La loi du 11 février 2005 fait de la LSF une langue à part entière et garantit aux parents de jeunes sourds une liberté de choix entre une communication bilingue - LSF et langue française - et une communication en langue française (avec ou sans langage parlé codé). La création des PASS (Pôles d'Accompagnement à la Scolarisation des élèves sourds) permet de scolariser des élèves sourds et malentendants en milieu ordinaire, quel que soit le mode de communication choisi par la famille. Ce choix s'exprime lors de l'élaboration par les parents du projet de vie de l'enfant et consécutivement à la demande déposée auprès de la Maison départementale des personnes handicapées (MDPH). Dans le cadre du PASS, la scolarisation des élèves sourds sur l'Académie d'Aix-Marseille, se faisait jusqu'à présent essentiellement dans le cadre d'un parcours « oralisant » - français oral et écrit - dans des structures d'accueil allant de la maternelle au collège. Concernant la scolarisation des élèves sourds dans le cadre d'un parcours bilingue - français écrit et LSF -, il n'existait pas de structures d'accueil en collège. Des structures étaient présentes dans des écoles maternelle et primaire (intégrations individuelles en classes ordinaires) et dans un lycée professionnel accueillant autant des élèves sourds oralisants que des élèves sourds bilingues. Les élèves bilingues étaient accueillis principalement dans un IRS (Institut Régional des enfants sourds). Depuis plusieurs années l'APES (Association des Parents des Enfants sourds) réclamait au niveau de l'Académie d'Aix-Marseille, l'ouverture de classes bilingues pour que leurs enfants et adolescents soient accueillis dans les écoles et les collèges. C'est dans ce contexte que s'est ouverte l'ULIS-Collège à Marseille.

[*] Aix-Marseille Université – France – teresa.dos-reis-assude@univ-amu.fr
[**] Aix-Marseille Université – France – karine.millon-faure@univ-amu.fr
[***] Aix-Marseille Université – France – jane.tambone@wanadoo.fr

La création de ce nouveau dispositif, le fait qu'il propose un parcours bilingue, a été pour nous une occasion de nous intéresser aux conditions de création et aux contraintes de la mise en place et du fonctionnement d'une nouvelle ULIS dans un collège.

Dans cette communication, nous nous intéressons à deux questions que nous formulons, dans un premier temps, de la manière suivante : Quels sont les obstacles que les élèves sourds rencontrent dans les apprentissages (notamment mathématiques) et sur quels leviers s'appuyer pour dépasser ces obstacles ? Nous allons aborder ces questions à partir d'un point de vue particulier qui est celui de la synchronisation du travail des élèves sourds avec celui des autres élèves. Pour cela, nous présenterons d'abord notre problématique et notre cadre théorique, ensuite nous préciserons la méthodologie avant de présenter quelques résultats.

## II. PROBLEMATIQUE ET CADRE THEORIQUE

La question des obstacles et des leviers pour les apprentissages d'élèves sourds a été abordée dans plusieurs travaux, notamment celui de Françoise Duquesne-Belfais (2007) sur lequel nous allons nous appuyer. Dans ce travail, l'un des obstacles mentionné est celui de la difficulté de conceptualisation des élèves sourds. Or Courtin (2002) considère que ces difficultés de conceptualisation ou d'abstraction ne sont pas dues à la surdité mais dérivent de l'insuffisance de situations expérimentées (notamment des situations sociales de communication et d'expression). Selon cet auteur, les élèves sourds ne sont pas suffisamment en contact avec une multiplicité de situations qui leur permette d'identifier des invariants pour pouvoir dégager des classes de situations propres au processus de conceptualisation. Par ailleurs, dans ce processus l'appui sur une langue semble essentielle, et la maîtrise de la LSF par des élèves sourds n'est pas toujours présente étant donné les différents choix des familles et le statut non reconnu à cette langue jusqu'à très récemment dans l'institution scolaire en France. Le choix d'une langue première comme la LSF (à l'instar d'autres langues) semble déterminant pour les élèves sourds comme moyen pour construire une représentation des situations, pour organiser, planifier une activité de résolution d'un problème, pour contrôler et réfléchir sur ses actions. Par contre, étant donné que la LSF ne possède pas un système d'écriture associé, les élèves apprennent le français écrit qui permet de mémoriser des résultats mais aussi d'organiser des raisonnements et de traiter des informations. Or des difficultés peuvent apparaître liés aux fonctionnements différents de la LSF en tant que langue gestuelle et visuelle, essentiellement iconique, et de la langue française écrite qui est une langue alphabétique, et linéaire. Les élèves sourds peuvent avoir des difficultés d'expression, une insuffisance de vocabulaire et une mauvaise maîtrise grammaticale (Duquesne-Belfais 2005, 2007).

L'un des leviers de l'enseignement des mathématiques dans un parcours bilingue est bien celui de s'appuyer sur la LSF, et de voir toutes les potentialités de cette langue dans la dialectique contextualisation/décontextualisation nécessaire au processus de conceptualisation. L'appui sur des situations contextualisées en nombre suffisant est important mais la décontextualisation est faite en utilisant des systèmes langagiers et non langagiers (Vergnaud 1990). La LSF peut être l'un de ces systèmes.

Selon Cuxac (2003), la LSF est une langue d'une iconicité importante : « Les séquences langagières sont en relation iconique très forte avec le point de départ extralinguistique de l'expérience réelle ou imaginaire à transmettre ». Cet auteur considère deux manières de dire dans cette langue : dire en donnant à voir, et dire sans montrer. La première a une visée d'illustration qui n'existe pas dans la deuxième. L'une des conditions est bien la présence des interlocuteurs pour pouvoir communiquer.

La scolarisation des élèves sourds dans des classes ordinaires avec l'aide d'un dispositif ULIS pose le problème de la communication entre les élèves sourds, le professeur et les autres élèves. Une insécurité des élèves sourds est celle de ne pas comprendre ce qui se passe autour d'eux et souvent les autres acteurs négligent l'importance d'expliquer ce qui se passe et pourquoi cela se passe (Duquesne-Belfais, Bertin 2005). A l'opposé, « l'enseignant dont les élèves sont sourds a souvent peur de l'implicite et de l'insécurité que porte en elle, par définition, une situation de recherche. Il peut alors être tenté, en LSF comme en français, de fournir trop d'explications qui ne laisseraient pas suffisamment à l'élève l'opportunité de faire des hypothèses et de mettre en œuvre un raisonnement par lui-même. » (Duquesne-Belfais, 2005, p.122). Comment faire en sorte que la communication puisse se faire en tenant compte des besoins des élèves sourds sans en faire trop ?

Dans le cas de notre dispositif, les élèves sont intégrés dans des classes ordinaires et les enseignants ne connaissent pas la LSF. La présence d'un interprète dans la classe ou dans le dispositif ULIS[1] est essentielle mais cela pose aussi des questions, par exemple celles des relations entre ce qu'on dit dans la classe et ce qui est traduit (dans les deux sens), celles du rapport aux savoirs enseignés des interprètes. Nous voulons interroger ce qui se passe en classe de mathématiques qui puisse être en relation avec cette contrainte communicationnelle liée à l'existence de deux langues.

Nos premières notions théoriques sont celles introduites par Chevallard (2010) à propos de l'organisation de l'étude autour d'un système didactique principal (SDP) et de plusieurs systèmes didactiques auxiliaires (SDA). La classe est le système didactique principal : le professeur et tous les élèves en font partie autour des enjeux de savoir. Pour scolariser les élèves sourds, un dispositif (ULIS) a été créé pour aider ces élèves relativement aux enjeux de savoirs. Le dispositif ULIS est un système didactique auxiliaire car il n'a pas d'enjeu de savoir propre et dépend de celui de la classe. C'est un dispositif didactique auxiliaire externe à la classe. Par contre, un autre système didactique auxiliaire existe à l'intérieur de la classe : celui constitué par les élèves sourds et l'interprète. Ce SDA interne à la classe est aussi dépendant du SDP et n'a pas d'enjeu de savoir propre.

Dans chaque système didactique (principal et auxiliaires), il existe diverses temporalités qui peuvent être propres ou communes à plusieurs systèmes. Le temps didactique (Chevallard et Mercier, 1987), défini comme le découpage d'un savoir dans une durée, est la norme de l'avancement des savoirs. Ce temps didactique pilote non seulement le travail dans la classe mais constitue une référence pour les systèmes didactiques auxiliaires. Plusieurs travaux ont montré qu'il n'y a pas d'avancement du temps didactique dans ces systèmes auxiliaires. D'autres temporalités existent dans ces systèmes (principal et auxiliaires), comme les temps d'apprentissage pour les élèves ou encore le temps praxéologique pour l'enseignement (Assude et alii, 2016). Le temps praxéologique est le temps qui correspond à l'avancement d'au moins l'une des composantes d'une praxéologie (Chevallard 1999), en se référant au quadruplet : « type de tâches, technique, technologie, théorie ».

---

[1] Les élèves suivent les enseignements dans leur classe avec les élèves entendants et ils disposent d'heures de regroupement dans le dispositif ULIS avec un enseignant spécialisé.

A partir de ces éléments théoriques, nous pouvons préciser notre question : quelles sont les relations entre les différents systèmes didactiques (principal et auxiliaires) du point de vue des temporalités ? Plus précisément, dans le cadre de la classe, quelles relations entre le SDP et le SDA interne à la classe ? Ou formulée autrement, comment sont synchronisés ces deux systèmes didactiques (SDP et SDA interne) ? Avant de donner des éléments de réponse à cette question, nous allons préciser notre cadre méthodologique ainsi que la situation proposée aux élèves en classe.

### III. ELEMENTS METHODOLOGIQUES RELATIFS A L'OBSERVATION

Pour étudier les relations entre le SDP et le SDA interne, nous avons observé une classe de mathématiques en 4$^{ème}$ qui accueille deux élèves sourds : Julien et Martin. Ces élèves sont sourds profonds ; ils suivent un parcours bilingue (LSF et français écrit), et bénéficient du dispositif ULIS avec quatre autres élèves sourds qui sont en sixième. Ces deux élèves ont quelques difficultés par rapport au français écrit, comme d'autres élèves sourds. Les moments de regroupement dans l'ULIS sont consacrés aux différentes disciplines, à la reprise ou finition de contrôles ou à l'aide aux devoirs.

Ces élèves suivent toutes les disciplines scolaires de 4$^{ème}$ sauf la deuxième langue vivante et la musique. La présence d'un interprète est constante dans toutes les disciplines et dans les moments de regroupement ULIS. Les quatre personnes qui assurent ces traductions n'ont pas toutes le même statut. Pour ce qui nous intéresse, l'interprète présent pour traduire dans la classe de mathématiques observée est le coordonnateur ULIS qui est un ancien professeur de mathématiques et qui connaît aussi la LSF.

Nous avons observé la classe de mathématiques de 4$^{ème}$ qui comporte 22 élèves pendant trois séances qui ont été filmées avec trois caméras : une caméra fixe qui filmait les deux élèves sourds, une caméra fixe qui filmait le tableau et une caméra mobile qui filmait le traducteur, ou l'enseignant lorsqu'il se déplaçait ou d'autres élèves. Nous avons aussi un entretien avec l'enseignant de mathématiques et un entretien avec l'interprète qui est aussi le coordonnateur ULIS[2].

Nous allons nous intéresser essentiellement à la première séance où les élèves devaient résoudre un problème de modélisation. La situation qui a été proposée à la classe est la suivante : *Combien de tours de pédale fait-on lorsqu'on parcourt une distance de 5km avec un vélo ?* Cette situation a été posée à l'oral, et trois phases ont été prévues par le professeur. La première phase est celle de la compréhension de l'énoncé et des contraintes du problème. La deuxième phase est celle d'un travail en groupe pour résoudre le problème. La troisième phase est celle de la mise en commun et de la correction du problème.

Ce problème pose quelques difficultés liées au processus de modélisation, notamment à la compréhension du fonctionnement d'un vélo. En voilà quelques-unes :

- Comprendre l'énoncé

---

[2] L'action du coordonnateur s'organise autour de 3 axes : l'enseignement aux élèves lors des temps de regroupement au sein de l'ULIS ; la coordination de l'ULIS et les relations avec les partenaires extérieurs ; le conseil à la communauté éducative en qualité de personne ressource. (Circulaire n° 2015-129 du 21-8-2015, BO n°31 du 27 août 2015)

- Dégager des variables importantes à partir de la situation réelle (relation entre tour de roue et périmètre du cercle, relation entre tour de pédale et tour de roue, vélo à 26 pouces, conversion pouce et mètres)
- Vocabulaire (Plateau/pignon/pédalier/vitesse)

L'enseignant a apporté un vélo en classe pour parler de tous ces éléments. Etant donné la difficulté de compréhension des énoncés des problèmes montrée dans d'autres travaux (Bonnet, Mangeret et Nowak, 2010), nous pouvons anticiper que les deux élèves sourds peuvent aussi se retrouver en difficulté pour entrer dans le milieu de la situation, même si c'est une situation contextualisée.

## IV. DESCRIPTION ET ANALYSE DE FAITS

Nous allons présenter des faits observés et des analyses du travail de la classe (SDP) en lien avec le SDA interne du point de vue de la synchronisation entre ces deux systèmes didactiques. Les faits observés concernent la compréhension de l'énoncé, les processus de modélisation et la maîtrise du vocabulaire indiqué auparavant.

### 1. *Une synchronisation globale…*

Durant les séances que nous avons observées, les élèves sourds et entendants se trouvent engagés en même temps dans les mêmes tâches : durant la première phase, tous ont tenté de s'approprier la situation et pour cela de comprendre le fonctionnement du vélo. Durant le travail de groupe, tous les élèves se sont lancés dans un travail de recherche. Enfin, tous les élèves ont pu participer à la phase de mise en commun.

En outre, nous avons pu observer que les élèves sourds communiquaient avec le reste de la classe. Ainsi même si le passage obligé par les traductions de l'interprète ralentit légèrement les échanges entre l'enseignant et les élèves, cela ne semble pas constituer une réelle entrave pour que les élèves sourds puissent intervenir. Martin n'hésite pas à lever le doigt pour répondre et l'enseignant l'interroge régulièrement. Julien, lui, intervient beaucoup moins souvent. En outre, lorsque le professeur se déplace entre les groupes pour vérifier que chacun avance dans la résolution du problème, il intervient aussi bien auprès des élèves sourds qu'entendants. Ainsi, même si elle nécessite l'intervention d'un interprète, la communication qui s'établit entre l'enseignant et les élèves sourds paraît semblable à celle que l'on peut observer avec des élèves ordinaires.

En ce qui concerne la communication entre élèves, les deux élèves sourds ne sont pas en dehors du reste de la classe. L'interprète souligne la bienveillance de leurs camarades à leur égard et la volonté de certains d'apprendre la langue des signes. Nous avons par ailleurs observé, lors de la deuxième séance, l'intervention d'une des élèves entendantes : ayant constaté que Martin n'avait pas ouvert le bon logiciel, elle a attiré son attention, grâce à des gestes, puis lui a montré son écran afin de lui faire comprendre son erreur. Ces échanges entre élèves sourds et élèves ordinaires semblent toutefois particulièrement rares. Ainsi, lorsque l'enseignant regroupe les élèves par 4 pour un travail de recherche, tous commencent à échanger entre eux, mis à part le groupe composé de Martin, Julien et deux élèves entendants. Immédiatement deux binômes se forment : les deux élèves sourds travaillent d'un côté, les deux élèves entendants de l'autre, sans qu'aucun échange entre eux ne se produise. L'interprète restera pourtant juste à côté d'eux durant tout le temps de la recherche, mais nul ne songera à le solliciter pour interagir avec l'autre binôme. De même il est à noter que l'échange observé en salle informatique, lors de la deuxième séance, se produit sans intervention de l'interprète, et l'on peut se demander si le fait de devoir recourir à une tierce personne pour communiquer entre eux ne constitue pas une difficulté pour les élèves. Toutefois, même s'il y a peu d'échanges entre eux, les élèves sourds et entendants semblent globalement synchrones.

## 2. ...mais quelques désynchronisations

Toutefois, en analysant nos séances de plus près, nous avons observé plusieurs épisodes où les élèves sourds prenaient du retard par rapport au le reste de la classe. Plusieurs causes peuvent être avancées pour expliquer ces désynchronisations.

Un premier facteur pourrait être lié au niveau de qualification en LSF de l'interprète. En effet ce dernier reconnaît qu'il n'est pas toujours à l'aise dans ses traductions et qu'il lui faut encore améliorer sa formation. Les contraintes de la traduction simultanée l'obligent à signer instantanément les expressions énoncées, ce qui lui pose parfois quelques difficultés. Ainsi, au début de la première séance, l'interprète ne sait comment traduire l'expression « tour de pédale » que l'enseignant vient d'utiliser. Il invente donc un geste qu'il espère suffisamment évocateur pour que les élèves le comprennent. Mais Martin et Julien s'imaginent qu'il vient d'utiliser le signe exprimant la marche ce qui complique la compréhension de la situation. L'éclaircissement de ce quiproquo nécessitera une mise au point entre les trois protagonistes afin de s'entendre sur la signification du geste de départ. De même un peu plus tard dans la séance, l'interprète cherche un signe pour désigner les vitesses du vélo (qu'il veut distinguer de son homophone « je roule à une vitesse de 5km/h ») mais le geste qu'il emploie n'est pas connu des élèves. Toutefois Julien comprend sa signification et en présente un autre, de sorte que les trois protagonistes parviennent tout de même à se comprendre. Certaines traductions de l'interprète se révèlent également ambiguës. Ainsi, dans la première séance, l'enseignant veut amener les élèves à trouver le nombre de tours de roue réalisés lorsqu'on fait un tour de pédale. Il explique donc qu'en un tour de pédale, « 42 dents passent dans la chaîne » à l'avant mais aussi à l'arrière, ce qui correspondra à 2 tours de pignon. Confronté à cette situation, l'interprète traduit la phrase « 42 dents passent dans la chaîne » par la phrase « il y a 42 dents à l'arrière ». Cette dernière phrase donne l'impression qu'il y a 42 dents sur le pignon alors qu'il y en a en fait que 21 (il faudra donc 2 tours de pédale pour que les 42 dents passent dans la chaîne), ce qui pourrait induire les élèves en erreur.

Une deuxième explication pourrait être liée au fait que la LSF est une langue visuelle et gestuelle : celle-ci ne peut se réaliser que si les deux interlocuteurs se voient et ont les mains libres. C'est ainsi que l'on remarque à plusieurs reprises l'interprète éprouver quelques difficultés pour attirer l'attention des élèves si ces derniers ne le regardent pas. Lors de la première phase, alors qu'il s'est approché du vélo pour leur montrer quelque chose, il doit faire de grands gestes avant d'entrer en communication avec eux. Par ailleurs, lorsque l'interprète traduit le discours de l'enseignant, les élèves sourds sont concentrés sur son visage et ses mains : ils ne peuvent donc regarder dans le même temps ce que l'enseignant de la classe est en train de tracer ou écrire au tableau. Enfin, durant les temps de copie, impossible pour les élèves sourds de profiter des éventuelles explications que l'interprète pourrait leur donner ou de lui poser des questions : ayant les yeux rivés sur le tableau ou sur leur cahier et les mains accaparées par le travail d'écriture, toute communication est impossible. Or, comme la plupart des enseignants, l'enseignant de la classe continue à parler pendant que les élèves copient et l'interprète ne peut donc pas traduire à Martin et Julien, ses propos.

Une troisième explication pourrait être celle d'un rapport difficile des élèves sourds au français écrit, comme le montrent d'autres études. Si les élèves sourds que nous observons savent a priori lire, leur rapport à l'écrit demeure plus compliqué que pour les élèves ordinaires. On peut en effet penser que pour ces derniers, la correspondance entre les phonèmes et les graphèmes facilite le déchiffrage. Les élèves sourds doivent associer à chaque mot écrit une idée sans pouvoir s'appuyer sur leur langue première qui est la langue des signes. Ceci pourrait expliquer les difficultés qu'ils rencontrent pour comprendre les mots qu'ils déchiffrent. En conséquence, même s'il a conscience que la lecture constitue pour eux un enjeu d'apprentissage essentiel, l'interprète doit régulièrement traduire en langue des signes une partie des supports proposés par les enseignants. En effet, sans cela, certains contre-sens graves pourraient empêcher leur entrée dans la tâche proposée. Ainsi, lorsque l'interprète montre le mot « pignon » écrit au tableau par l'enseignant, et demande « vous comprenez ce qu'il y a écrit ? », Martin répond sans hésiter « oui, c'est un oiseau ». En effet, cet élève a associé les mots « pignon » et « pigeon » qui présentent une écriture assez proche. Cette difficulté posée par la compréhension du français écrit amène l'interprète à ne pas se contenter de traduire les échanges entendus dans la classe. En effet, les élèves sourds doivent apprendre non seulement le geste-signe associé à une notion donnée mais également l'écriture du mot correspondant, sans qu'aucun lien évident n'existe entre ces deux signifiants. L'interprète tente donc d'accompagner ces apprentissages en écrivant régulièrement sur une ardoise les mots correspondants aux notions nouvelles qu'il vient de signer. Il signale aussi la polysémie des mots ou l'existence d'homographes. Ainsi lorsque l'enseignant parle des dents de l'engrenage, l'interprète commence par signer cette notion, puis il fait remarquer aux élèves que ce terme s'écrit comme les dents situées à l'intérieur de la bouche.

### 3. *Resynchronisation : des bulles de compréhension*

Nous voyons donc que les éléments pouvant engendrer des retards du SDA par rapport au SDP sont nombreux. Toutefois, comme la classe avance, le SDA interne doit aussi avancer et réussir à se resynchroniser par différents moyens. Nous avons ainsi pu constater dans les séances observées que l'interprète ne traduisait pas toujours l'intégralité des propos qui s'échangeaient dans la classe : il lui est arrivé de ne sélectionner que les informations les plus importantes pour la résolution du problème posé en occultant le reste.

Un autre dispositif peut permettre de resynchroniser les élèves sourds avec le reste de la classe. En effet, durant la séance, à plusieurs reprises, des discussions, entre les deux élèves sourds, ou entre les élèves sourds et leur interprète, s'instaurent indépendamment des échanges qui se produisent dans le reste de la classe. Ce sont des temps où l'interprète va pouvoir s'assurer de la compréhension de certaines notions évoquées et enrichir les explications données par l'enseignant. Martin et Julien discuteront également de leurs interprétations respectives de la situation proposée. Ainsi, pendant que le reste de la classe échange sur le fonctionnement du vélo dans la première séance, Martin et Julien parlent de la conversion cm/pouces puis du calcul du périmètre. On a alors l'impression que Martin, Julien et l'interprète sont comme coupés du reste de la classe dans une sorte de sas isolé leur permettant d'éclaircir certains points incontournables pour accéder à l'activité mathématique attendue. C'est pourquoi nous parlerons de *bulles de compréhension*. Toutefois, durant ces apartés, le déroulement de la séance se poursuit et certains échanges entre les élèves ordinaires et l'enseignant ne seront pas traduits. L'interprète se doit donc d'estimer l'intérêt de ces deux sources d'informations afin de décider si à un moment donné, l'isolement des élèves sourds constitue un appui susceptible de favoriser leur compréhension et de permettre de les resynchroniser avec le reste de la classe ou au contraire un handicap en raison des informations perdues.

V.  CONCLUSION

Les séances que nous avons observées se composent d'une succession de moments de synchronisation et de désynchronisation : à certains moments les élèves sourds sont parfaitement en phase avec le reste de la classe, à d'autres ils sont en décalage, soit à cause de ralentissements dus à la langue des signes ou à leur rapport à l'écrit, soit grâce à l'apparition de bulles de compréhension qui peuvent permettre de resynchroniser ces élèves avec les autres. Ces cycles dé-resynchronisations n'empêchent pas les élèves sourds d'être engagés dans la résolution du problème, de vouloir comprendre de quoi il s'agit, d'essayer de trouver des techniques qui permettent de trouver la réponse au problème. On peut donc se demander si ces désynchronisations constituent de réels obstacles pour les apprentissages des élèves sourds : finalement, a-t-on besoin d'être toujours synchrone ?

Notre travail montre plutôt des réponses nuancées à ces questions. Certes tout dépend du moment où ces désynchronisations ont lieu, et à propos de quoi. Dans notre cas, la première phase qui permet de comprendre les contraintes du problème notamment celles concernant le fonctionnement du vélo est essentielle pour la dévolution du problème aux élèves. Le temps didactique avance car il y a des objets nouveaux qui sont essentiels pour s'engager dans la recherche d'une technique. Les resynchronisations sont alors nécessaires, et le rôle de l'interprète est important pour condenser une discussion, pour l'apport d'une information. Nous avons ainsi observé que l'interprète revenait sur les éléments écrits au tableau et montrait certaines parties du vrai vélo pour que les élèves puissent comprendre la relation entre tour de pédale et tour de roue. Le temps praxéologique a avancé puisque les élèves sourds, tout comme les élèves entendants, ont pu s'engager dans la résolution du problème. Pendant le travail en groupe, le SDA interne d'une part et les autres élèves d'autre part ont fonctionné d'une manière autonome. Différents temps praxéologiques co-existaient dans la classe puisque les techniques des groupes se mettaient en place. A la troisième étape, la mise en commun a permis d'institutionnaliser les techniques qui permettaient d'aboutir à la solution du problème. Pendant ce temps, le SDA interne était synchronisé avec la classe. Le temps didactique a ainsi avancé pour la classe mais aussi pour le SDA interne.

Ce travail nous permet de donner une réponse aux obstacles et aux leviers sur lesquels s'appuyer pour favoriser les apprentissages des élèves sourds. La présence d'un SDA interne à la classe peut être un obstacle si les désynchronisations associées à la présence de deux langues empêchent d'avoir les informations essentielles pour s'engager dans le travail mathématique. Cela peut se faire de manières très diverses, par des bulles de compréhension, par des cycles de désynchronisation/resynchronisation plus ou moins rapides. Le fait que dans la séance observée l'interprète était auparavant un enseignant de mathématiques a sûrement facilité la sélection des informations présentées aux élèves. Par la suite, nous voulons étudier ce qui se passe si l'interprète n'est pas spécialiste de la discipline : est-ce que dans ce cas le temps didactique avance en classe et aussi dans le SDA interne ?